\documentclass[11pt]{amsart}
\usepackage{amsbsy,amssymb,amscd,amsfonts,latexsym,amstext,delarray,
amsmath,graphicx}
\input xypic

\newtheorem{thm}{Theorem}[section]
\newtheorem{prop}[thm]{Proposition}

\newtheorem{defn}[thm]{Definition}
\newtheorem{rem}[thm]{Remark}

\numberwithin{equation}{section}

\def\C{{\mathbb C}}

\renewcommand{\P}{{\mathbb P}}

\def\Z{{\mathbb Z}}
\def\R{{\mathbb R}}

\def\cL{{\mathcal L}}

\def\cP{{\mathcal P}}

\newcommand{\ie}{{\it i.e.,\/}\ }
\newcommand{\eg}{{\it e.g.\/}\ }
\newcommand{\cf}{{\it cf.\/}\ }

\def\text{\hbox}

\def\Hom{{\rm Hom}}

\title{Open string theory and planar algebras}
\author[Ceyhan]{\"Ozg\"ur Ceyhan}
\author[Marcolli]{Matilde Marcolli}
\address{\"O.~Ceyhan: Max--Planck Institut f\"ur Mathematik  \\
Vivatsgasse 7 \\ Bonn, D-53111 Germany} 
\email{ceyhan@mpim-bonn.mpg.de}
\address{M.~Marcolli: Mathematics Department, California Institute of Technology \\
1200 E.California Blvd.f \\ Pasadena, CA 91125, USA} 
\email{matilde@caltech.edu}

\begin{document}

\maketitle

\begin{abstract}
In this note we show that abstract planar algebras are
algebras over the topological operad of moduli spaces of 
stable maps with Lagrangian boundary conditions, which 
in the case of the projective line are described in 
terms of real rational functions. These moduli spaces
appear naturally in the formulation of open string theory
on the projective line. We also show two geometric ways to
obtain planar algebras from real algebraic geometry, one
based on string topology and one on Gromov-Witten theory.
In particular, through the well known relation between
planar algebras and subfactors, these results establish 
a connection between open string theory, real algebraic
geometry, and subfactors of von Neumann algebras. 
\end{abstract}

\section{Introduction}

The purpose of this paper is to show that planar algebras
arise as algebras over an operad of moduli spaces of stable
maps to $\P^1$ with Lagrangian boundary conditions, which
can be described in terms of real algebraic curves. 
In particular, the results presented in this paper can be 
interpreted as a connection between open string theory on 
$\P^1$ and Jones' theory of subfactors of von Neumann
algebras, using as intermediate steps the relation between open string 
theory on $\P^1$ and certain moduli spaces $R_g(\P^1, d)$ of maps to 
$\P^1$ with Lagrangian boundary conditions, combined with the 
main result we prove here, which relates the latter 
to the theory of planar algebras developed in \cite{Jones}. 
The connection to subfactors can then be seen by invoking the 
result of \cite{Popa} and its reformulation in terms of planar 
algebras of \cite{GJS} and \cite{KoSu2}. 

The main point involved in our description of planar algebras
in terms of real algebraic geometry is an identification 
of (weighted) planar tangles as the combinatorial datum that 
encodes the components of the moduli spaces $R_g(\P^1, d)$ of stable 
$M$-maps to $\P^1$ with Lagrangian boundary conditions. The
terminology $M$-maps here refers to the fact that they are realized 
by maximal real algebraic curves. We also describe the 
compositions of planar tangles and the trace map in this geometric
setting. More precisely, we prove the following statement.

\begin{thm}\label{planRg}
Abstract planar algebras are topological operads of $R_g(\P^1,d)$.
\end{thm}

In particular, by restricting to the case with $g=0$, we obtain the
case of the Temperley-Lieb algebras.

The first two sections of this paper respectively review the relation
between open string theory on $\P^1$ and the moduli spaces
$R_g(\P^1,d)$, as well as the well known relation between subfactors and
planar algebras. The main original contribution of this paper is in 
\S \ref{mainSec} and \S \ref{AlgSec}. In \S \ref{mainSec} we connect 
the two previous topics by proving Theorem \ref{planRg} and
in \S \ref{AlgSec} we then describe two different methods,
both based on real algebraic geometry, for constructing 
planar algebras as representations of the planar operad, 
the first based on algebraic loop spaces and string 
topology and the other based on a real version of 
Gromov--Witten theory.

%%%%%%%%%%%%%%%%%%%%%%%%%%%%%%%%%%%%%%%%%%%%%%%%%%%%%%%%%%%%%%%%
%  
%                moduli spaces again =)
%
%%%%%%%%%%%%%%%%%%%%%%%%%%%%%%%%%%%%%%%%%%%%%%%%%%%%%%%%%%%%%%%%
%\newpage

\section{Real algebraic geometry and open string theory}

\subsection{Complex and real curves, and bordered Riemann surfaces }

Recall that a bordered Riemann surface of type $(g,h)$ is a Riemann surface 
with boundary, with $g$ handles and $h$ boundary components, oriented 
according to the orientation induced by the complex structure on the Riemann 
surface. We also denote, as in \cite{AlGr}, \cite{KatzLiu} the complex double 
of a bordered Riemann surface $\Sigma$ by $\Sigma_\C$, with $\sigma$ the 
antiholomorphic involution on $\Sigma_\C$. 

\subsubsection{Maximal real curves}

A maximal real curve, denoted $M$-curve, is a real algebraic curve
$(\Sigma_\C,\sigma)$ with the maximal number of connected components of
the real part $\Sigma_\R$ of $\sigma$. By Harnack's bound this maximal number is
$g + 1$ for a curve of genus $g$. 

The real part $\Sigma_\R$ of a real structure $\sigma$ divides
$\Sigma_\C$ into two 2-dimensional discs, respectively denoted by
$\Sigma^+$  and $\Sigma^-$, minus a set of interior (open) discs
$D_1,\ldots,D_g$, having $\Sigma_\R$ as their common boundary in
$\Sigma_\C$. 
The real structure $\sigma$ interchanges $\Sigma^\pm$, and the complex
orientations of $\Sigma^\pm$ induce two opposite orientations on $\Sigma_\R$, 
called its complex orientations.  The quotient 
$\overline{\Sigma} := \Sigma_\C/\sigma$ is isomorphic to $\Sigma^+$. 
Here $\Sigma^+\simeq \overline{\Sigma}$ is a bordered Riemann surface 
and the real algebraic curve $(\Sigma_\C, \sigma)$ is its complex double, 
\cf \cite{Sepp}.

%------------------moduli of real maps--------------------------

\subsection{Open strings on $\P^1$ and the moduli spaces of real maps}

We present here briefly the setting of open string theory on $\P^1_\C$ 
with a Lagrangian submanifold $\P^1_\R$. We
especially focus on the role of the moduli space of stable maps with
Lagrangian boundary conditions.

\subsection{M-maps and their moduli spaces}

We define stable $M$-maps following a setting similar to that of \cite{KatzLiu}.

\begin{defn}\label{stableMmap}
An  $M$-map (of degree $d$) is a pair $(\Sigma^+,f^+)$ 
consisting of the following data. 

\begin{itemize}
\item A bordered Riemann surface $\Sigma^+$ whose complex double 
$(\Sigma_\C,\sigma)$ is an $M$-curve of genus $g$, and with 
boundary components $\partial \Sigma = \partial D_* \cup \partial 
D_1 \cup \cdots \cup \partial D_g$. 

\item A map $f^+$ of degree $d$ which is the restriction of a real 
stable morphism $f: \Sigma_\C \to \P^1_\C$, such that the singularities 
of $f$ are all of degree two.

\item A set consisting of a critical point of $f$ for each non-empty 
component of  $\Sigma_\R$.  
\end{itemize}
An  $M$-map is  stable if it does not admit any infinitesimal automorphisms.
\end{defn}

The moduli space  $R_g(\P^1,d)$ is the space of isomorphism classes 
of stable $M$-maps. The geometric and topological properties of the 
space $R_g(\P^1,d)$ has been extensively studied in \cite{Liu}. The
construction in \cite{Liu} is in fact more general; it considers an
arbitrary symplectic manifold $X$ and its Lagrangian manifold $L$ as
the target of stable maps.

%Notice that the moduli space $R_g(\P^1,d)$ is non-compact.

%------------------planar tangles--------------------------

\section{Planar algebras and subfactors}

\subsection{Planar algebras}
We recall here briefly the basic definitions and facts of the theory
of planar algebras developed in \cite{Jones}. We follow closely the 
short survey given in \cite{Bisch}.

\subsubsection{Planar tangles}\label{SecTangle}
A planar $k$-tangle $T$ consists of the unit disc $D\subset \C$ 
together with a finite collection of discs $D_1 , D_2 , \ldots, D_g$
inside $D$. The boundaries of $D$ and of each interior disc $D_i$ are
decorated by an even number of marked points, $2k$ points on $\partial
D$ and $2k_i$ points on each $\partial D_i$. The interior of
$D\smallsetminus \bigcup_i D_i$ also contains a collection of
non-intersecting strings, which are either closed strings or have as
boundary the marked points on the $\partial D \cup \partial D_i$. Each
marked point lies on the boundary of one of these strings.  

The complementary region $D \smallsetminus (\{\text{strings}\} \bigcup_i
D_i)$ admits a checkerboard black and white coloring as well as 
a choice of a white region at each $D_i$.

\begin{center}
\begin{figure}
\includegraphics[scale=0.4]{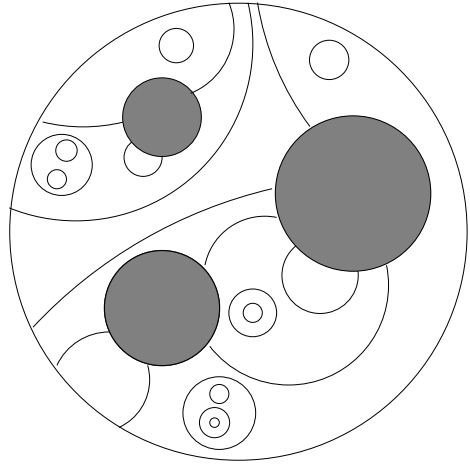}
\caption{Planar tangle \label{FigPlanDiag}}
\end{figure}
\end{center}

\subsubsection{Planar operad}
Two planar tangles $T$ and $S$ can be composed whenever the
number of marked points on the boundary of $S$ matches the number of
marked points on the boundary of one of the interior discs $D_j$ of
$T$. The composition $T \circ_j S$ is then given by gluing a rescaled
copy of $S$ in place of the interior of the disc $D_j$, so as to match
the shadings and the marked white regions. 
This operation is well-defined, with the coloring and choice of white region
eliminating any possible rotational ambiguity. Also the result only
depends on the isotopy classes. 

\begin{center}
\begin{figure}
\includegraphics[scale=0.4]{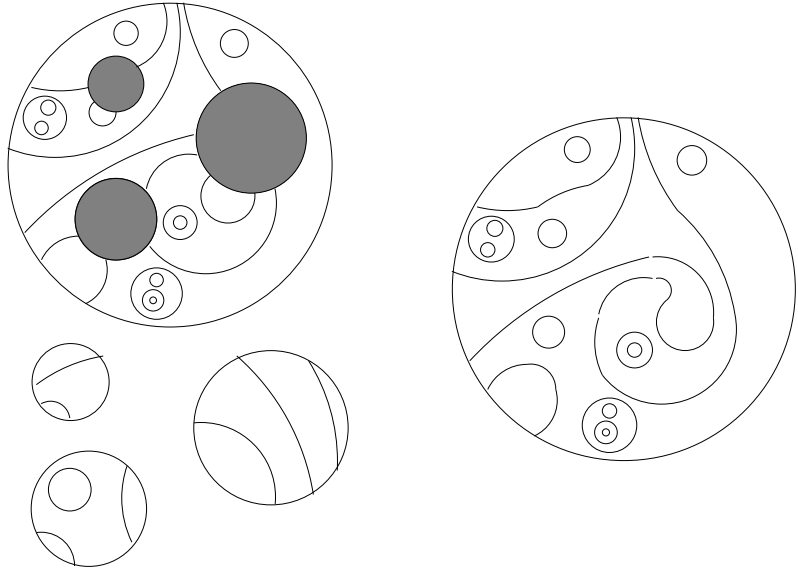}
\caption{Planar operad: compositions of tangles \label{FigPlanOpd}}
\end{figure}
\end{center}

The planar operad $\cP$ consists of all orientation-preserving 
diffeomorphism classes of planar k-tangles that fix the boundary
$\partial D$. The structure of operad is given by the compositions of
tangles, defined as above. An example of compositions is given in
Figure \ref{FigPlanOpd}. In the figures \ref{FigPlanDiag} and
\ref{FigPlanOpd} the black/white coloring and 
the marked white regions are not shown for simplicity.

\subsubsection{Planar algebras}
One then defines planar algebras (\cite{Jones}, \cite{Bisch}) as
algebras over the planar operad $\cP$, in the sense described in
\cite{May}. 

This means (see \eg \cite{Bisch}) that a planar algebra $\cP$ is a 
family of vector spaces $\{ V_k \}_{k>0}$ together with a morphism
$Z$ from the planar operad $\cP$ to the (colored) operad $\Hom$ of
multilinear maps between vector spaces. 

\subsubsection{Partition function and trace} The morphism $Z : \cP \to
\Hom$ from the planar operad to the operad of multilinear maps of
vector spaces has the following properties.

Given a tangle $T_\ell$, where $2\ell$ is the number of points on
$\partial D$, one obtains a multilinear map 
\begin{equation}\label{ZTell}
Z(T_\ell): \otimes_{i=1}^n V_{\ell_i} \to V_\ell, 
\end{equation}
satisfying the composition property
\begin{equation}\label{Zcomp}
Z(T \circ_j S) =Z(T)\circ_j Z(S).
\end{equation}

Moreover, one has the following properties, which give $Z$ an
interpretation as a ``partition function'' associated to a planar
algebra. 

\begin{center}
\begin{figure}
\includegraphics[scale=0.4]{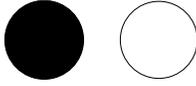}
\caption{The planar tangles $T_{0,\emptyset,b/w}$ \label{FigT0bw}}
\end{figure}
\end{center}

$\bullet$ Normalization: let $T_{0,\emptyset,b/w}$ be the colored
discs of Figure \ref{FigT0bw}. Then
\begin{equation}\label{Znorm}
Z(T_{0,\emptyset,b/w})=1.
\end{equation}

\begin{center}
\begin{figure}
\includegraphics[scale=0.4]{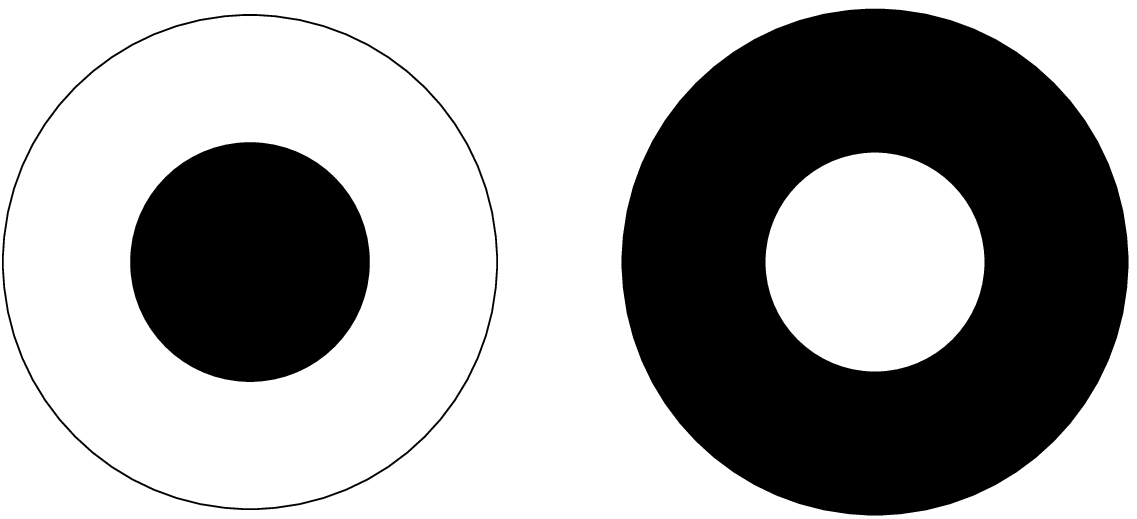}
\caption{The planar tangles $T_{0,\emptyset,\cL,b/w}$ \label{FigTLbw}}
\end{figure}
\end{center}

$\bullet$ Parameters: Consider the case of planar tangles as in Figure
\ref{FigTLbw}, given by a simple curve with no boundary inside $D$ and
the two possible choices of coloring. Then the corresponding linear
maps under $Z$ are of the form
\begin{equation}\label{Zdelta}
Z(T_{0,\emptyset,\cL,b})=\delta_1,  \ \ \ \
Z(T_{0,\emptyset,\cL,w})=\delta_2, 
\end{equation}
for two parameters $\delta_i$. 

\begin{center}
\begin{figure}
\includegraphics[scale=0.5]{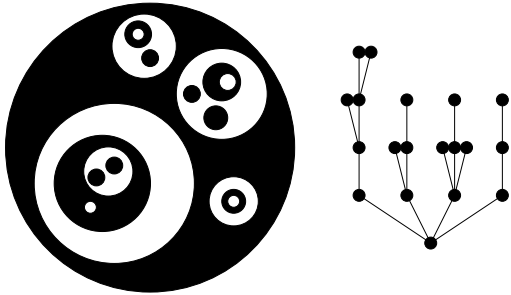}
\caption{Ovals in planar tangles encoded by rooted trees \label{FigOvals}}
\end{figure}
\end{center}

These two rules give a procedure to {\em eliminate ovals} from a
planar tangle, replacing them in the image under $Z$ by a
multiplicative factor of the form $\delta_1^m \delta_2^n$,
depending on the parameters $\delta_i$. This is illustrated in one
example in Figure \ref{FigOvals}, where one encodes a configuration
of ovals inside a planar tangle via a rooted tree and computes the 
corresponding multiplicative factor accordingly. 

%------------------standart invariants--------------------------

\subsection{Subfactors and planar algebras}

Our main focus in this paper is on a geometric framework for planar
algebras based on real algebraic geometry and open string theory.
However, we mention briefly the well known connection
between planar algebras and the theory of subfactors of von Neumann 
algebras (\cite{Bisch},\cite{GJS}, \cite{Jones},\cite{KoSu1}, 
\cite{KoSu2}, \cite{Popa}),
since, in light of our interpretation of planar algebras it leads 
naturally to possible constructions of subfactors from data of
moduli spaces in real algebraic geometry.

In \cite{Popa}, Popa showed that one can associate to a subfactor 
a {\em standard invariant}, which is a planar algebra. Not all
planar algebras arise as the standard invariant of a subfactor,
but one can characterize those that have this property, as in
\S 4 of \cite{Jones}.

A subfactor planar algebra is a planar algebra for which all
the vector spaces $V_k$ are finite dimensional, with $\dim V_{0,b/w}=1$,
and with $\delta_1=\delta_2 \neq 0$, which has an involution on each $V_k$
induced by the reflection of tangles with the property that the partition
function $Z$ is a sesquilinear form with respect to this involution. 
The reconstruction theorem, in the form given in \cite{GJS}, shows
how to associate to a subfactor planar algebra with $\delta>1$ 
a subfactor $M_0\subset M_1$ obtained from a tower of type II$_1$ factors  
$M_k= Gr_k(V)=\oplus_{n\geq k} V_k$, with a faithful tracial state $tr_k$,
so that the standard invariant of the subfactor, constructed as in \cite{Popa} 
by considering the relative commutants $M_0^\prime \cap M_k$, gives
canonical identifications $M_0^\prime \cap M_k \simeq V_k$ which induce
a morphism of involutive planar algebras.

%\newpage
%------------------moduli of real maps--------------------------

\section{Planar algebras and open string theory}\label{mainSec}

%------------------moduli space & tangles-----------------------

\subsection{The moduli space $R_g(\P^1,d)$ and planar tangles}

We first define  weighted planar tangles, which are
planar tangles with additional decorations, and we use 
them to distinguish the connected components of the 
moduli space $R_g(\P^1,d)$.

\begin{defn}\label{weightTangle}
Let $T$ be a planar tangle as in \S \ref{SecTangle}. Then, 
$T$ is called weighted if each interior string and each boundary 
segment is decorated  with a non-negative integer weight. A 
weighted tangle is of weight $d$ if the sum of weights for all
strings and boundary segments is $d$.
\end{defn}

The following statement is a slight generalization of a similar 
result in \cite{NSV} for genus zero covers of $\P^1$.

\begin{prop}\label{CompRgTangles}
There is a 1-1 correspondence between the set of connected components 
of the moduli space $R_g(\P^1,d)$ and the set of weighted planar tangles.
\end{prop}

\proof 
We associate a tangle to each stable $M$-map $(\Sigma^+,f^+) \in 
R_g(\P^1,d)$ by simply pulling back the real line 
$\P^1_\R= \R \cup \infty$ by the map $f^+$,
\begin{eqnarray}
(\Sigma^+,f^+) \mapsto T := (f^+)^{-1} (\P^1_\R).
\end{eqnarray}

Let $(\Sigma_\C, f)$ be the complex double of $(\Sigma^+,f^+)$.  
Let $Crit(f)$ be the set of critical values of $f$, and let 
$Crit_\R (f)$ and $Crit_\C (f)$ denote respectively the set of 
real critical values and its complement. To determine weights,
consider an arbitrary point in $x \in \P^1_\R \setminus Crit_\R (f)$.
For each string (or boundary segment), we define the function 
$w(x)$ as the number of preimages of $x$ lying on this string 
(or boundary segment). Then the weight of this string (or boundary 
segment) is the minimum of $w(x)$. Note that the weight of a closed
string in $\Sigma^+$ coincides with the multiplicity of $f^+$
restricted to this string.

We first check that this indeed gives a weighted planar tangle.
Since $\Sigma^+$ is half of a maximal real curve of genus $g$,
it is topologically a disc minus a collection of $g$ discs.
The preimage $(f^+)^{-1} (\P^1_\R)$ gives the strings inside
$\Sigma^+$ or segments in the boundary of $\Sigma^+$.

We now check the condition that these tangles are expected to satisfy.
These strings intersect only in the boundary 
$\partial \Sigma^+$ of $\Sigma ^+$. We have the following possibilities.

\vspace{.2cm}
\paragraph{\it Strings inside $\Sigma^+$}
If the strings had intersections in $\Sigma^+ \setminus 
\partial \Sigma^+$, then such an intersection point $z$ would be a critical 
point with real critical value. Hence, the same is also true for the 
conjugate $\sigma(z) \in \Sigma_\C$ in the complex 
double $\Sigma_\C$ of $\Sigma^+$. However, this contradicts 
the genericity condition of $f: \Sigma \to \P^1_\C$ 
(see Defn.~\ref{stableMmap}).

\vspace{.2cm}
\paragraph{\it Strings and the boundary $\partial \Sigma^+$}
All critical points of $f$ which have real critical values 
must be real. If $z$ is a critical point with a real critical 
value, then $\bar z$ is a critical point with the same critical 
value. Therefore, $z = \bar z$, since $f$ is generic. For each 
such critical point, $f^{-1}(\P^1_\R)$ contains exactly four 
arcs in $\Sigma_\C$ incident to it. Two of these arcs are the 
arcs lying in $\partial D_i \subset \R\Sigma$, while the other 
two interchange under the involution $\sigma$. In particular, 
the other endpoints of these two arcs coincide.

Next, we need to show that the number of 
critical points at each boundary component $\partial D_i$ of 
$\Sigma^+$ is even, as required for the data to define a planar 
tangle. Note that the set of all critical points of the real map 
$f: \Sigma \to \P^1_\C$ consists of an even number of points. The 
order of this set is in fact  $2(d+g-1)$ due to the  
Riemann--Hurwitz formula. On the other hand, the real maps 
$(\Sigma_\C,f)$ degenerate when their critical points collide, 
in particular when a complex conjugate pair of critical points 
degenerates to a real point of $\Sigma$. Under such a degeneration, 
it can only split into a pair of real ramification points lying 
in the same real component $\partial D_i$. Since there exist 
stable maps $f: \Sigma_\C \to \P^1_\C$ without any real 
ramification points, we show that the number of critical points 
is even for each boundary component $\partial D_i$ via a simple 
induction on these types of degenerations.

We need to show that, for the tangle we associated to elements 
of a connected component of $R_g(\P^1,d)$, the 
the isotopy class does not change. It is clear that the isotopy class
of these tangles will not change under small perturbations and
it can only change through intersections of strings. As we have 
already observed above, any intersection of strings violates the 
genericity condition. Namely, changes of isotopy class of tangles
can happen when one passes through the discriminant locus to another
component of $R_g(\P^1,d)$.

Finally, to show that we have a bijection, we then need to
show that any arbitary  tangle of weight $d$ can be obtained in 
this way. For a given such tangle, an element of $R_g(P^1,d)$
can be constructed by gluing a set of weighted pairs of pants
which is in fact $\Sigma^+ \setminus \{\textrm{strings}\}$. The 
gluing prescription of weighted pants is given in
 Thm 1 in \cite{NSV} based on \cite{Nat}.
\endproof

In the following, we denote by $C_T$ the connected component of 
$R_g(\P^1,d)$ corresponding to a given weighted tangle $T$.

%------------------composition of tangles--------------------------

\subsection{Sewing stable $M$-maps}
\label{sec_sewing}

Let $(\Sigma^+_1,f_1^+)$  and $(\Sigma^+_2,f_2^+)$ be a pair of 
stable $M$-maps such that the restriction of $f_1^+$ to the boundary
component $\partial D_{i_1} \subset \Sigma^+_1$ agrees with the 
restriction of $f_2^+$ onto $\partial D_{i_2} \subset \Sigma^+_2$. 
Assume that the boundary components $\partial D_{i_1}$ and $\partial D_{i_2}$ 
carry opposite orientations. Then we can sew the $M$-maps
$(\Sigma^+_k,f_k^+)$ along their boundaries $\partial D_{i_k}$ 
for $k=1,2$.

Let $R_{g_1}(\P^1,d_1) \times_{ij} R_{g_2}(\P^1,d_2)$ denote the 
space of pairs described above. This space is in fact a fiber product
in the following way.

Recall that the space of {\em algebraic loops} is the
space of morphisms $f: \P^1_\R \to \P^1_\R$.
The moduli spaces $R_{g_k}(\P^1,d_k), k=1,2$ admit evaluation maps
to $L^{alg}(\P^1_\R)$,
\begin{eqnarray}
ev_{i_k}: R_{g_k}(\P^1,d_k) \to L^{alg}(\P^1_\R), \ \ 
(\Sigma^+_k,f_k^+)  \mapsto f^+_k |_{\partial D_{i_k}}  \nonumber
\end{eqnarray}
Then the space  described above is the fibered product 
\[ R_{g_1}(\P^1,d_1) \times_{L^{alg}(\P^1_\R)} R_{g_2}(\P^1,d_2). \]

The sewing operation described above provides us with a morphism 
\begin{equation}\label{sewcompose}
R_{g_1}(\P^1,d_1) \times_{ij} R_{g_2}(\P^1,d_2) \to 
R_{g_1+g_2-1}(\P^1,d_1+d_2).
\end{equation}

The induced map on the sets of connected components of these
spaces defines a corresponding map at the level of weighted tangles, which 
agrees with the composition of tangles.

%------------------sewign boundary--------------------------

\subsubsection{Sewing boundary segments}
There is a similar sewing along the boundary segments which fits better
with the example that we discuss in \S \ref{TLsec} below.

Let $(\Sigma^+_1,f_1^+)$  and $(\Sigma^+_2,f_2^+)$ be a pair of 
stable $M$-maps such that the restriction of $f_1^+$ to the boundary
segment $I_{j_1} \subset \partial D_{i_1} \subset \Sigma^+_1$ agrees 
with the restriction of $f_2^+$ onto $I_{j_2} \subset \partial D_{i_2} 
\subset \Sigma^+_2$. If  $\partial D_{i_1}$ and $\partial D_{i_2}$ 
carry opposite orientations, then  we can sew the $M$-maps
$(\Sigma^+_k,f_k^+)$ along their boundary segments $\partial I_{j_k}$ 
for $k=1,2$ and obtain a new $M$-map of degree $d_1 + d_2$.

The space of such pairs is also a fiber product. The appropriate
subspaces of the moduli spaces $R_{g_k}(\P^1,d_k), k=1,2$ admit 
evaluation maps to the space of algebraic paths $P^{alg}(\R)$, 
which is the space of morphisms $f: \R \to \R$. Note that 
this map is not defined for all components $C_T$ of $R_{g_k}(\P^1,d_k)$
since the tangles $T$ in general need not have boundary segments
(\eg if there is no critical point on a boundary component of $\Sigma^+$,
then it can only map into the algebraic loop space as above).
Then, by using the appropriate subspaces of $R_{g_k}(\P^1,d_k)$, we 
obtain the space of such pairs as a fiber product similar to the above 
construction (\ref{sewcompose}).

%------------------fianlly, the trace--------------------------

\subsection{The trace}
\label{sec_trace}

Let $(\Sigma^+,f^+)$ be in $C_T \subset R_g(\P^1,d)$ and let
$S$ be the closed string which is the common boundary of a pair of weighted 
pants $P_1, P_2$  in $\Sigma^+$. Let $\{\partial D_1,\dots, \partial D_l,S\}$ 
and $\{\partial D_{l+1},\dots, \partial D_k,S\}$  be the sets of  boundaries 
of $P_1$, and $P_2$ respectively. Let the $w_i$ be the weights of $\partial D_i$
and $w_s$ be the weight of $S$.

We first note that the deformations of such a pair of weighted pants and their
maps to $\P^1$ are given by the fiber product
\begin{eqnarray}
R_{l}(\P^1, w_1 +\cdots + w_l) \times_s R_{k-l}(\P^1, w_{l+1} +\cdots + w_k)
\nonumber
\end{eqnarray}
which is determined by the evaluation map $(P_i, f_i) \mapsto f |_{S}$. As we
have already seen above, there is a morphism of this product into the
moduli space $R_{k-1} (P^1, w_1 +\cdots + w_k)$. This map in fact removes the
string $S$ in $\Sigma^+$ \ie it provides the trace operator in geometric
terms in the above setting. This map is given by the gluing of a pair of weighted pants $P_1$ and $P_2$ as in~\cite{NSV}.

%------------------and ,the invlution--------------------------

\subsection{The involution}
\label{sec_inv}
We can define an involution $*$ by replacing the complex structure
$J$ of $\Sigma^+$ by $-J$. This is equivalent to replacing $\Sigma^+$
with the other half $\sigma (\Sigma^+)$ of $\Sigma_\C$. This operation 
reverses the orientations of the strings and checkerboard shadings 
of the corresponding planar tangle.

\subsection{Planar algebras and the moduli space of real maps $R_g(\P^1,d)$}
The collection of connected components of $R_g(\P^1,d)$ provides
a topological operad $\cP_{\P^1}$. The following statement then follows
directly from the discussion of the previous subsections.

\begin{thm}
The operad of (weighted) planar tangles is the topological operad $\cP_{\P^1}$.
\end{thm}

%------------------temperley lieb--------------------------

\subsection{An example: Temperley-Lieb algebras}\label{TLsec}
\begin{center}
\begin{figure}
\includegraphics[scale=.5]{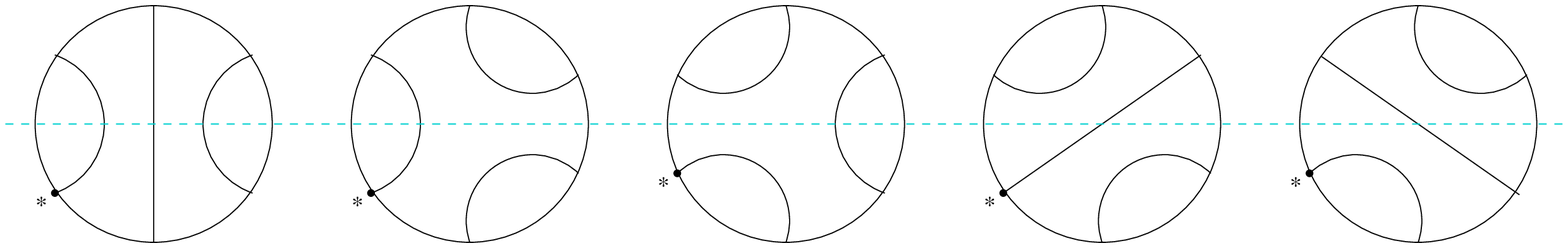}
\caption{The bases of  $TL_3$ \label{tl3}}
\end{figure}
\end{center}

Consider the moduli space $R_0(\P^1,d)$ of $M$-maps of genus zero. 
Such $M$-maps have $2d-2$ critical points. Then, the tangles that
distinguish the components $C_T$ of $R_0(\P^1,d)$ are the planar
tangles on a disc with $2k$ points on their boundaries where 
$0 \leq 2k \leq 2d-2$. In other words, these tangles connect the
first $k$ points starting at the marked critical point in $\Sigma^+$
to the second $k$ points without having any crossings. These are in fact 
the tangles that generate the Temperley-Lieb algebras (modulo
the trace). For instance, if we consider the  case with $d=4$, we obtain
the tangles that generate $TL_1, TL_2$ and $TL_3$; see Figure
\ref{tl3} for the bases of $TL_3$.

\begin{center}
\begin{figure}
\includegraphics[scale=.5]{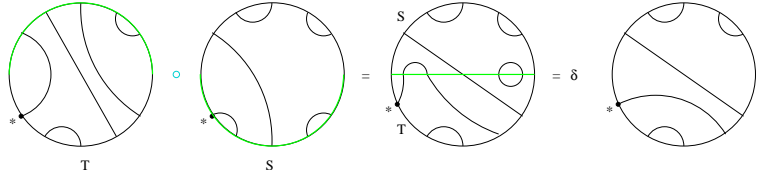}
\caption{The product of tangles $S$ and $T$ \label{tlproduct}}
\end{figure}
\end{center}

The product of two such tangles is obtained by sewing the boundary 
segments between the $k$th and the $2k$th point of the first tangle and the $1$st 
and the $k$th point of the second tangle. Figure \ref{tlproduct} illustrates
an example of a such a composition of tangles.

%------------------genus one--------------------------

\begin{rem}{\rm 
There is an equivalent description of Temperley-Lieb algebras
by using the moduli space of $M$-maps of genus one. In this case,
one needs to consider the tangles in a cylinder $\Sigma^+$, and the
composition is given by sewing the boundary components. However, 
the numbers of critical points at each boundary of the cylinder
need not be equal. Therefore, either one has to consider a more 
general algebra (which contains Temperley-Lieb algebra as a subalgebra) 
or one needs to restrict oneself to a subspace of $R_1(\P^1,d)$ that 
contains only those components $C_T$ that correspond to the tangles 
of the Temperley-Lieb algebra.}
\end{rem}

%------------------Algebras over operad--------------------------
 
\section{Algebras over planar operad}\label{AlgSec}
In this section, we discuss two possible geometric constructions of 
algebras over the planar operad $\cP$. The first one follows closely 
the approach of Cohen and Godin's in contrucing string topology operations 
\cite{CohenGodin}. The second one is in the spirit of Gromov-Witten theory 
\cite{KontManin}, in the real algebraic geometry setting (see for 
instance \cite{Ce}).

%------------------string topology--------------------------

\subsection{Planar algebras in string topology}
\label{sec_str_top}
As we have already observed in \S \ref{sec_sewing}, the moduli space 
$R_g(\P^1,d)$ admits evaluation maps 
$$ ev_i: R_g(\P^1,d) \to L^{alg}(\P^1_R), \ \ \  i=1,\dots,g+1. $$
By using these evaluations, we can pull back classes from 
$H^*(L^{alg}(\P^1_R))$. Note that, like the usual topological loop space $L S^1$,
the algebraic loop space $L^{alg}(\P^1_R)$ has infinitely many connected 
components, each determined by the degree of the loops. However, in contrast to
the topological case, the components of $L^{alg}(\P^1_R)$ are finite dimensional 
and are not contractible. Therefore they carry nontrivial cohomology classes.

Let $\gamma_1,\dots, \gamma_{g+1} \in H^*_c(L^{alg}(\P^1_R))$. Then, we
define correlators
\begin{eqnarray}
\langle \gamma_1,\dots,\gamma_{g+1} \rangle_T := ev^* (\gamma_1) \wedge 
\cdots \wedge ev^* (\gamma_{g+1})  
\end{eqnarray}
which take values in the cohomology with compact support  $H_c^* (C_T)$
of the component $C_T \subset R_g(\P^1,d)$ corresponding to the tangle
$T$.

Let $\{e_a\}$ be a basis for $H^*_c(L^{alg}(\P^1_R))$ and let $g^{ab}$ be the
intersection form.  We  use the
above invariants of the loop space $L^{alg}(\P^1_R)$ to give a representation 
of the operad $\cP$ of planar tangles in the following way. 
The image of $T$ under the morphims $Z: \cP \to \Hom$ is  
\begin{equation}
\label{eqn_g_product}
T \mapsto Z(T) := \left(\gamma_1 \otimes \dots \otimes \gamma_{g}  \mapsto 
\sum_{a,b}  \langle \gamma_1,\dots,\gamma_{g}, e_a \rangle_T \  g^{ab}  \ e_b \right).
\end{equation}
The composition of these higher products are evident from the definitions
of these products in (\ref{eqn_g_product}) and from the compositions of tangles
described in \S \ref{sec_sewing}.

In this construction, the trace can be given by using the procedure described in \S 
\ref{sec_trace}. Let $T$ be a planar tangle and let $S$ be closed in $T$. 
In such a case, we think of $T$ as a composition of two tangles $T_1$ and $T_2$ 
as in \S \ref{sec_trace}.
Then, we obtain this composition as in equation (\ref{Zcomp}).
In particular, if the closed string bounds a disc (rather than a more 
complicated weighted pant), then it plays the role of trace, and it is
given as in \eqref{Zdelta}.

%------------------string topology--------------------------

\subsection{Planar algebras via Gromov-Witten theory}
\label{sec_gw}
An alternative way to obtain representations of $\cP$ is to use
ideas from Gromov-Witten theory of $\P^1$. Let $T$ be a tangle having at
least one marked point at each boundary disc. In this case, the 
evalution map 
\begin{eqnarray}
ev: C_T \to Conf(\P^1_\R,g+1) \times \R Conf(\P^1,k),\ \
\text{where} \ \ k=n-g+1, \nonumber
\end{eqnarray}
maps $(\Sigma^+,f^+)$ to its critical values in $\P^1$. The image
of the evalution map is the space of $\Z/ 2\Z$-equivariant distinct 
(unordered) point configurations in $\P^1$ which we denote by 
$\R Conf(\P^1,k)$. Moreover, if we restrict ourselves to the 
specified critical points at each boundary component, the evalution 
map takes its values in the configuration space of ordered point
configurations $Conf(\P^1_\R,g+1)$. 

Let $\alpha_1,\dots, \alpha_k \in H^*_c(\R Conf(\P^1,k))$ and 
$$ \gamma_1,\dots, \gamma_{g+1} \in H^*(Conf(\P^1_\R,g+1)).$$ 
We then define correlators by setting
\begin{eqnarray}
\langle \langle \alpha_1,\dots, \alpha_k; \gamma_1,\dots, \gamma_{g+1}
\rangle \rangle_T := 
ev^* (\alpha_1 \wedge \cdots \wedge  \alpha_k \wedge
\gamma_1 \wedge \cdots \wedge \gamma_{g+1}).  \nonumber
\end{eqnarray}
These take values in the cohomology with compact support  $H_c^* (C_T)$
of the component $C_T \subset R_g(\P^1,d)$ corresponding to the tangle
$T$. Then, by using the same idea as in (\ref{eqn_g_product}), we define  
the higher products corresponding to the tangles $T$, which depend on
the classes $\alpha_i$, by setting
\begin{eqnarray}
\label{eqn_gw_product} \nonumber
T \mapsto Z(T) := \left(\gamma_1 \otimes \dots \otimes \gamma_{g}  \mapsto 
\sum_{a,b}  \langle \langle \alpha_1,\dots, \alpha_k; \gamma_1,\dots, 
\gamma_{g}, e_a \rangle \rangle_T \  \eta^{ab}  \ e_b \right).
\end{eqnarray}
Here $\eta^{ab}$ is the intersection matrix of $H^*(Conf(\P^1_\R,g+1))$.
The composition of these higher products and the trace operators are  
evident and are given as in the previous case described in \S \ref{sec_str_top}

\begin{rem}{\rm 
The requirement on the existence of marked points may seem artificial, 
but at present we do not have a convenient substitute for it.
However, the above setting adapts well to the examples of Temperley--Lieb
algebras and Fuss--Catalan algebras (see \cite{Bisch}).}
\end{rem}

\begin{rem}{\rm 
Both situations discussed above, based on string topology or on Gromov-Witten
theory, can be further enriched by using the tautological classes of the moduli
space $R_g(\P^1,d)$. This might lead to additional interesting examples.}
\end{rem}

\bigskip

{\bf Acknowledgement.} The first author would like to thank 
Vaughan Jones for useful discussions. 
The second author is partially supported by NSF-grants 
DMS-0651925 and DMS-0901221. Part of this work was carried
out during the authors stay at the Max Planck Institute for
Mathematics, which we thank for the hospitality and for
support. 

%------------------references--------------------------
%------------------------------------------------------


\begin{thebibliography}{99}

\bibitem{AlGr} N.L.~Alling, N.~Greenleaf, {\em Foundations of the
theory of Klein surfaces}, Springer, 1971. 

%\bibitem{At} M.~Atiyah, {\it Topological quantum field theories}. 
%Inst. Hautes Etudes Sci. Publ. Math. Vol.68 (1988), 175--186.

\bibitem{Bisch} D.~Bisch, {\em Subfactors and planar algebras}, ICM
2002 Beijing, Vol.III, 1--3.

\bibitem{Ce} O.~Ceyhan, {\em Towards quantum cohomology of real varieties}, 
arXiv:0710.0922.

\bibitem{CohenGodin} R.~Cohen,  V.~Godin, 
{\em A polarized view of string topology.}  
Topology, geometry and quantum field theory,  
127--154, London Math. Soc. Lecture  
Note Ser., 308, Cambridge Univ. Press, Cambridge, 2004.

%\bibitem{Cost} K.~Costello, {\em The $A_\infty$ operad and the moduli
%space of curves}, arXiv:math/0402015.

%\bibitem{EG} A.~Eremenko, A.~Gabrielov, {\em Rational functions with 
%real critical points and the B. and M. Shapiro conjecture in real 
%enumerative geometry}.  Ann. of Math. (2)  155  (2002),  no. 1, 105--129. 

\bibitem{GJS} A.~Guionnet, V.F.R.~Jones, D.~Shlyakhtenko, {\em Random
matrices, free probability, planar algebras and subfactors},
arXiv:0712.2904. 

\bibitem{Jones} V.F.R.~Jones, {\em Planar algebras, I},
math.QA/9909027.

\bibitem{KatzLiu} S.~Katz, C.C.M.~Liu, {\em Enumerative geometry of
stable maps with Lagrangian boundary conditions and multiple covers of
the disc},  Adv. Theor. Math. Phys.  5  (2001),  no. 1, 1--49. 

\bibitem{KoSu1} V.~Kodiyalam, V.S.~Sunder, {\em A complete set of
numerical invariants for a subfactor}, J. Funct. Anal. 212 (2004)
1--27. 

\bibitem{KoSu2} V.~Kodiyalam, V.S.~Sunder, {\em From subfactor planar
algebras to subfactors}, arXiv:0807.3704.

\bibitem{KontManin} M. Kontsevich, Y.U.~Manin, {\em Gromov-Witten classes, 
quantum cohomology, and enumerative geometry.}  Comm. 
Math. Phys.  164  (1994),  no. 3, 525--562.

\bibitem{Liu} C.C.M.~Liu, {\em Moduli of J-holomorphic curves with
Lagrangian boundary conditions and open Gromov-Witten invariants for
an $S^1$-equivariant pair}, arXiv:math/0210257. 

\bibitem{May} J.P.~May, {\em Definitions: operad, algebras and
modules}, Contemporary Mathematics, 202 (1997) 1--7.

\bibitem{NSV} S.~Natanzon, B.~Shapiro, A.~Vainshtein, {\em 
Topological classification of generic real rational functions}. 
J. Knot Theory Ramifications 11 (2002), no. 7, 1063--1075. 

\bibitem{Nat} S.~Natanzon, {\em Topology of 2--dimensional coverings 
and meromorphic functions on complex and real algebraic curves}. 
Selecta Math. Sovietica 12 (1993),  251--291. 

\bibitem{Popa} S.~Popa, {\em 
An axiomatization of the lattice of higher relative commutants 
of a subfactor}. Invent. Math. 120 (1995), no. 3, 427--445. 

\bibitem{Sepp} M.~Sepp\"al\"a, {\em Moduli spaces of stable real
algebraic curves}, Ann. Sci. \'Ecole Norm. Sup. (4) 24 (1991)
519--544. 

\end{thebibliography}
\end{document}